\documentstyle[11pt]{amsart}
\title[A finiteness theorem $\dots$]{
A finiteness theorem for  low-codimensional \\ nonsingular subvarieties 
of quadrics}
\author{Mark Andrea A.  de Cataldo}

\address{Department of Mathematics,
Washington University in St. Louis,
Campus Box 1146,
St. Louis, Missouri 63130-4899}
 \thanks{$^*$ L. Fania and G. Ottaviani have recently informed 
the author that they have settled the
case $n=6$ in the positive.}

\email{mde@@math.wustl.edu}


\newtheorem{tm}{Theorem}[section]
\newtheorem{lm}[tm]{Lemma}
\newtheorem{pr}[tm]{Proposition}
\newtheorem{rmk}[tm]{Remark}
\newtheorem{cor}[tm]{Corollary}
\newtheorem{ex}[tm]{Example}

\font\tenmsb=msbm10
\font\sevenmsb=msbm7
\font\fivemsb=msbm5

\newfam\msbfam
\textfont\msbfam=\tenmsb
\scriptfont\msbfam=\sevenmsb
\scriptscriptfont\msbfam=\fivemsb
\def\Bbb#1{{\fam\msbfam #1}}

\font\teneufm=eufm10
\font\seveneufm=eufm7
\font\fiveeufm=eufm5
\newfam\eufmfam
\textfont\eufmfam=\teneufm
\scriptfont\eufmfam=\seveneufm
\scriptscriptfont\eufmfam=\fiveeufm

\newcommand\ci{\cite}
\newcommand\s{\sigma}

\newcommand\zed{{\Bbb Z}}
\newcommand\pn[1]{{\Bbb P}^{#1}}

\newcommand{\Q}[1]{{\cal Q}^{#1}}

\newcommand\blacksquare{{\hspace*{\fill} $\Box$}} 
\newcommand\odix[1]{ {\cal O}_{#1} }
\newcommand\odixl[2]{ {\cal O}_{#1}({#2}) }

\begin{document}
\centerline{\bf TO APPEAR IN TRANS. A.M.S.}
\bigskip
\
\bigskip

\maketitle

\begin{abstract}
We prove that there are only finitely many families of codimension two 
nonsingular
subvarieties of quadrics $\Q{n}$ which are not of general type,
 for $n=5$ and $n\geq 7$.
We prove a similar statement also for the case of higher codimension.
\end{abstract}

\section{Introduction}
\label{boundednessofthreefoldsnotgentype}

 There are only finitely many
families of codimension two nonsingular subvarieties not of general type
 of the projective spaces $\pn{n}$,
for $n\geq 4$; see \ci{e-p} and \ci{boss1}. More generally,
a similar statement holds for the case of higher codimension; see
\ci{sc}.

In this paper we concentrate on the case of codimension two
subvarieties of quadrics.
Our  main result  is Theorem \ref{sigh}: there are only finitely
many families of nonsingular codimension two subvarieties
 not of general type  in the quadrics $\Q{n}$, $n=4,$ $5$ or $n\geq 7$.
The case  $n=4$ is proved in 
\ci{a-s}, \S6. The case of $n=5$ is at the heart of the paper;
the main tools are the semipositivity of the normal bundles of
nonsingular subvarieties of quadrics, the Double Point Formula,
 the Generalized Hodge Index Theorem, bounds
for the genus of curves
on $\Q{3}$,  Proposition \ref{roth} and
 Corollary \ref{finiteonhyp}.
The case $n=6$ 
is still open$^*$.
The case of codimension two with $n\geq 7$ is covered by
 Theorem \ref{boundednessonqn},
 which hinges upon the result of
\ci{sc}; it also gives  a finiteness result in codimension 
bigger than two
in the same spirit as \ci{sc}.

\smallskip
The paper is organized as follows.
Section \ref{prelim} records, for the reader's convenience,  
some results used
in the paper.  A generalization of a  lifting criterion of Roth's
is contained in section
\ref{rothtypelifting}; we shall need the particular case
expressed by Proposition \ref{rothpr}.
Section  \ref{boundednessforn>=7} deals with higher codimensional
cases.
Sections \ref{chTHFD} and \ref{boundednessonq5} 
 are modeled  on   \ci{boss1}. 
Section \ref{chTHFD}  contains the lengthy proof of 
Theorem \ref{polybound} and of its Corollary
\ref{finiteonhyp}.
Section \ref{boundednessonq5}
contains the proof of Theorem \ref{sigh}.

\smallskip \noindent
{\bf Notation and conventions.}
Our basic reference is \cite{ha}.
We work over any algebraically closed field of characteristic zero.
A quadric $\Q{n}$,
here, is a nonsingular hypersurface of degree two in the  
projective space
$\pn{n+1}$.
Little or no distinction is made between line bundles, 
associated sheaves of sections
and Cartier divisors; moreover the additive and tensor notation
are both used.

\smallskip
\noindent
{\bf Acknowledgments.} It is a pleasure to thank our Ph.D. advisor
A.J. Sommese, who has suggested to us to study $3$-folds in $\Q{5}$.

\section{Preliminary material}
\label{prelim}

\begin{pr}
\label{easybound}
{\rm (Cf. \ci{decastharris} or \ci{thesis}; for the case
of $d> 2k(k-1)$ see \ci{a-s}, \S6.)}
Let $C$ be an integral curve of degree $d$ contained in an 
integral surface
of degree $2k$ in $\Q{3}$. Then the following bound holds
for the genus $g$ of $C$:
$$
g-1\leq \frac{d^2}{4k} + \frac{1}{2} (k-3)d. 
$$
\end{pr}

\begin{pr}
\label{boundasep}
{\rm (Cf. \ci{a-s}, Proposition $6.4$.)}
Let $C$ be an integral curve in $\Q{3}$, not contained in any
surface in $\Q{3}$ of degree strictly less than $2k$. Then:
$$
g-1\leq \frac{d^2}{2k} +\frac{1}{2}(k-4)d.
$$
\end{pr}

\noindent

Let $S$ be a nonsingular surface in $\Q{4}$, $N$ its normal bundle, 
  $C$ a nonsingular hyperplane
section of $S$, $g$ its genus, $d$ its degree.
Let $V_s \in |{\cal I}_{S,\Q{4}}(s)|$, where $s$ is some
 positive integer, 
be an integral hypersurface and 
$\mu_l:=c_2(N(-l))=$ $(1/2)d^2 +l(l -3)d-2l(g-1)$, 
$\forall l\in \zed$.

\begin{lm}
\label{epas}
{\rm (Cf. \cite{thesis}, Lemma $2.35$.)}
In the above situation:
$$
0\leq \mu_{s}\leq s^2 d.
$$

\end{lm}

The following proposition follows immediately from Theorem \ref{roth}
when the ambient space ${\cal P}^{n+2}$ is chosen to be a 
quadric $\Q{n+2}$. 

\begin{pr}
\label{rothpr}
Let $X$ be an integral subscheme of degree $d$ and codimension 
two in $\Q{n}$,
 $n\geq 4$. Assume that for the general hyperplane section
$Y$ of $X$ we have
$$
h^0(\Q{n-1}, {\cal I}_{Y,\Q{n-1}}(\s))\not= 0,
$$
for some positive integer $s$ such that $d>2{\s}^2$.
Then
$$
h^0(\Q{n}, {\cal I}_{X,\Q{n}}({\s}))\not= 0.
$$
\end{pr}

\medskip
Let $X$ be a degree $d$,  nonsingular $3$-fold in $\Q{5}$,  $L\simeq 
{\odixl{\pn{6}}{1}}_{|X}$, $S$ be the   surface
general hyperplane section of $X$, $C$ the general curve section of $S$
and $g$ the genus of $C$. 
 Using the
Double Point Formula (cf. \ci{fu}) for the embedding
$X\hookrightarrow \Q{5}$, we  get
the  following formul\ae\   for
$K_X\cdot L^2$, $K_X^2\cdot L$, $K_X^3$  as functions
of $d$, $g$, $\chi ({\cal O}_X)$, $\chi ({\cal O}_S)$:

\begin{equation}
\label{KL2}
K_X\cdot L^2=2(g-1)-2d,
\end{equation}
\begin{equation}
\label{K2L}
K_X^2\cdot L=\frac{1}{4}d^2 +\frac{3}{2}d -8(g-1) +6\chi ({\cal O}_S),
\end{equation}
\begin{equation}
\label{K3}
K_X^3=-\frac{9}{4}d^2+\frac{27}{2}d+gd +18(g-1)-30\chi 
({\cal O}_S)-24\chi ({\cal O}_X),
\end{equation}

Finally we record the expression for the Hilbert polynomial of $X$,

\begin{eqnarray}
\label{chioxt}
\chi  ({\cal O}_X (t))  =   \frac{1}{6}dt^3+[\frac{1}{2}d
- \frac{1}{2}(g-1)]t^2+    
  [\frac{1}{3}d - \frac{1}{2}(g-1) +\chi ({\cal O}_S)]t +
\chi ({\cal O}_X).   
\end{eqnarray}

For the details concerning the above formul\ae \  see \ci{thesis}, \S1.

\subsection{A Roth-type lifting criterion}
\label{rothtypelifting}

 If the general curve section 
of a degree $d$ linearly normal surface $S$ in $\pn{4}$ lies on a
surface of degree $\s$ in $\pn{3}$, then $S$ lies on 
some hypersurface
of degree $\s$, provided $d>\s^2$ (cf. \cite{ro}).
A generalization of this fact to codimension two integral linearly 
normal subschemes
of $\pn{n}$, $n\geq 4$ has been known for some time.

In this section we generalize Roth's lifting criterion 
to a larger class of spaces;
see Theorem \ref{roth} and Example \ref{examples}. 
The proof does not require the concept of linear normality, which was
virtually automatic in the case that Roth considered.

The proof given below was inspired by  
\cite{a-s}, Lemma $6.1$.

\noindent
First we fix some notation.
Let  ${\cal P}^{n+2}$ be a nonsingular projective variety of dimension
$(n+2)$, $n\geq 2$,  $L=\odixl{{\cal P}^{n+2}}{1}$ an ample and spanned 
line bundle on it with
$\delta:= L^{n+2}$. Assume that 
$Pic({\cal P}^{n+2})\simeq \zed [L]$.
 Let $X^n$ be an integral subscheme of
${\cal P}^{n+2}$ of dimension $n$ and
$d:=L^{n}\cdot X$.
Denote by ${\cal P}^{i+2}$   the intersection of 
$(n-i)$ general elements of $|L|$ and  by $X^i$ the intersection of
the same elements  of $|L|$ with $X^n$. 

\begin{tm}
\label{roth} Assume that the natural restriction maps below 
are surjective
$\forall m$:
$$
 \rho_m:= H^0({\cal P}^{n+2},mL)\to H^0({\cal P}^{n+1},
 mL_{|{\cal P}^{n+1}}).
$$
If  
$h^0({\cal I}_{X^{n-1},{\cal P}^{n+1}}(s))\neq 0$ for some
$s$ such that
$d>\delta s^2$,
then $h^0($ ${\cal I}_{X^{n},{\cal P}^{n+2}}(s))$ $\neq 0$.

\noindent
If $s$ is the minimum such number then
 $h^0({\cal I}_{X^{n},{\cal P}^{n+2}}(s))=1$.

\end{tm}

\noindent
{\em Proof.}
Let us assume that we have proved the theorem for
$$s=\sigma:=\min
\{t\in {\Bbb N}|\quad  h^0 ( {\cal I}_{X^{n-1},{\cal P}^{n+1}} (t) )
\neq 0
\};
$$
we call $\s$ the {\em postulation} of $X^{n-1}$. 
Then the theorem
holds also $\forall s\geq \sigma$. We can thus assume, without
loss of generality, that $s=\sigma$.

\smallskip
\noindent
Pick any $V_{\sigma} \in |{\cal I}_{X^{n-1}, {\cal P}^{n+1}}(\s)|$.

\medskip
\noindent
CLAIM. {\em $V_{\sigma}$ is integral}. 
This follows easily
from the minimality
of $\sigma$ and the fact that, under our assumptions,
 $Pic({\cal P}^{n+1})\simeq
\zed [L_{|{\cal P}^{n+1}}]$.

\medskip
\noindent
CLAIM. 
{\em 
$V_{\sigma}$ 
is  the unique element of 
$|{\cal I}_{X^{n-1},{\cal P}^{n+1}}(\sigma)|$}.
 By contradiction, assume that we have two distinct
$V^i_{\sigma}$. 
By the above claim they are both integral. By an easy
 Bertini-type argument we see that intersecting everything with
$n$ general members of $|L|$ we get two distinct integral curves
$W^i_{\sigma}\in |\odixl{{\cal P}^2}{\sigma}|$ containing
$X^0=\{d$ points$\}$. Since the curves do not have common components
we see that
$d\leq W^1_{\sigma}\cdot W^2_{\sigma}=\delta {\sigma}^2$; 
the intersection product here
is on ${\cal P}^2$. 
This  is a contradiction and the claim is proved.

\smallskip
\noindent
Let us choose a general line 
$\ell \subseteq$ 
$|L|^{\vee}$. Define
$\tilde{\cal P}$ 
to be the blowing up of 
${\cal P}^{n+2}$ along the intersection
of all the members of $\ell$. Denote by
$p$ and $q$
 the natural projections
to $\ell$ and ${\cal P}^{n+2}$, respectively. By intersecting 
with general 
elements of 
$|L|$ we get the following diagram, where $Y^i$ denotes
$q^{-1}X^i$:
\medskip

$
\begin{array}{lllllllll}
 \hspace{1cm} Y^0 &
 \subseteq & 
Y^1 & 
\subseteq  &
 \ldots  &
\subseteq &
 Y^n & 
\subseteq &
 \tilde{\cal P}\stackrel{p}{\longrightarrow} \ell  \\
\hspace{1cm} \downarrow & \ & \downarrow & \ & \  & \  
& \downarrow & \ & \downarrow q \\
  \hspace{1cm} X^0            & \subseteq & X^1 & \subseteq & 
\ldots & \subseteq &  X^n
& \subseteq & {\cal P}^{n+2}.              
\end{array}
$
\medskip

\noindent
We  have the following injections, where, for simplicity
(and by abuse) of notation,
we denote a twist by
$q^*\odixl{{\cal P}^{n+2}}{\s}$ simply by
a twist by ${\s}$:
$$
{\cal I}_{Y^n,\tilde{\cal P}}(\sigma)\to
{\cal I}_{Y^{n-1},\tilde{\cal P}}(\sigma)\to
\ldots
\to
{\cal I}_{Y^1,\tilde{\cal P}}(\sigma)\to
{\cal I}_{Y^0,\tilde{\cal P}}(\sigma),
$$
so that, applying $p_*$, we obtain the following injections:
$$
p_*{\cal I}_{Y^n,\tilde{\cal P}}(\sigma)\to
p_*{\cal I}_{Y^{n-1},\tilde{\cal P}}(\sigma)\to
\ldots
\to
p_*{\cal I}_{Y^1,\tilde{\cal P}}(\sigma)\to
p_*{\cal I}_{Y^0,\tilde{\cal P}}(\sigma).
$$
The existence of $V_{\sigma}$, for a general point of $\ell$, ensures
that
$p_*{\cal I}_{Y^n,\tilde{\cal P}}$ $(\sigma)$ is not the zero sheaf.
Since $p$ is dominant and the ideal sheaves
${\cal I}_{Y^i,\tilde{\cal P}}$ are torsion free, we see that 
the sheaves
$p_*{\cal I}_{Y^i,\tilde{\cal P}}(\sigma)$
are torsion free $\forall i$. But $\ell$ is a smooth curve, so that 
these sheaves are actually locally free. The uniqueness statement,
which was shown above, implies  that these sheaves are actually 
line bundles
on $\ell$. Since each of the above injections has torsion free 
cokernel
on $\ell$ we deduce that they all are isomorphisms, i.e.:
$$
p_*{\cal I}_{Y^n,\tilde{\cal P}}(\sigma)\simeq
p_*{\cal I}_{Y^{n-1},\tilde{\cal P}}(\sigma)\simeq
\ldots
\simeq
p_*{\cal I}_{Y^1,\tilde{\cal P}}(\sigma)\simeq
p_*{\cal I}_{Y^0,\tilde{\cal P}}(\sigma)\simeq
\odixl{\ell}{\tau},
$$
for some $\tau \in \zed$.

\noindent
By contradiction, assume $\tau<0$. Then
%\begin{eqnarray*}%
%0 & = & h^0(p_* {\cal I}_{Y^{n-1},\tilde{\cal P}}(\sigma))=%
%h^0({\cal I}_{Y^{n-1},\tilde{\cal P}}(\sigma))=         \\%
 %& = &h^0(q_*{\cal I}_{Y^{n-1},\tilde{\cal P}}(\sigma))=%
%h^0({\cal I}_{X^{n-1},{\cal P}^{n+2}}(\sigma)).%
%\end{eqnarray*}%
$$
0  =  h^0(p_* {\cal I}_{Y^{n-1},\tilde{\cal P}}(\sigma))=
h^0({\cal I}_{Y^{n-1},\tilde{\cal P}}(\sigma))=         
   h^0(q_*{\cal I}_{Y^{n-1},\tilde{\cal P}}(\sigma))=
h^0({\cal I}_{X^{n-1},{\cal P}^{n+2}}(\sigma)).
$$

\noindent
By our assumptions we can lift a section defining
$V_{\sigma}$ to a non zero element of 
$H^0($ ${\cal I}_{ X^{n-1},{\cal P}^{n+2} }$ $(\sigma))$.  
This contradiction proves $\tau \geq 0$. 

\noindent This proves the first assertion
 of the theorem. As to the second we need to prove that
$\tau =0$. 
 $\tau $ being strictly positive  would violate the usual uniqueness.
\blacksquare

\begin{rmk} 
{\rm We used the surjectivity of the restriction maps only for 
$m=\sigma$.}
\end{rmk}

\begin{rmk}
{\rm The cases
$({\cal P}^{n+2},L)\cong $ $ (\pn{n+2},$ $\odixl{\pn{n+2}}{1})$,
$ (\Q{4},$
$\odixl{\Q{4}}{1})$
seem to be well known. See for example \cite{a-s}, \cite{m-r}, 
and of course
\cite{ro}. However in the case of projective space it seems that the
linear normality of $X$ was usually  required; after Zak's 
theorem on tangencies
linear normality is automatic, for a nonsingular $X$, unless 
$ n=4$ and
$X$ is the Veronese surface.}
\end{rmk}

\begin{ex}
\label{examples}
{\rm The variety ${\cal P}^{n+2}$ can be, for example, a
 projective space,
a nonsingular
complete intersection or  a Grassmannian; in all these cases
 $L$ is the hyperplane bundle
for the natural embedding.
But it can also be chosen to be
  a Fano variety,
 of index $r$, with $-K_{\cal P}=rL$, $L$ generated by global 
sections
 and $Pic({\cal P})\simeq \zed$ (this is always the case if
$r>n/2$), some weighted complete intersections or, more generally, low 
degree branched coverings of projective
spaces \cite{la1} or Grassmannians \cite{ki}. 
In the last batch of examples, $L$ does not need to be very ample.}
\end{ex}

The following gives a lifting criterion in any codimension;
see  \cite{m-r}. Again linear normality is not required.

\begin{cor} Let $X^{\nu}$ be an integral subscheme of  
${\cal P}^{n+2}$ of 
dimension
$\nu$, $X^{\nu -1}$ the intersection of  $X$ with a general member of
$|L|$, $\s$ the postulation of $X^{\nu -1}$.
Assume that
$h^0({\cal I}_{X^{\nu -1}}(\s))=1$ and that $\rho_{\s}$ is surjective.
Then $h^0({\cal I}_{X^{\nu}}(\s))=1$.
 \end{cor}

\section{Finiteness on $\Q{n}$, $n\geq 7$}
 \label{boundednessforn>=7}

In this section we remark that, for a nonsingular variety of dimension
$\nu \geq \frac{n+3}{2} $ in $\Q{n}$ not of general type, the bound 
$d\leq 2 n^{n - \nu}$ holds. This gives the finiteness of 
the corresponding
 number of families. 
  
We thank M. Schneider for
pointing out to us that the result of this section could 
be proved along the lines of
his paper
\ci{sc}.

\begin{tm}
\label{boundednessonqn}
There are only
finitely many components of the Hilbert scheme
of $\Q{n}$ corresponding to  nonsingular subvarieties of dimension 
$\nu \geq \frac{n+3}{2} $ which are not of general type.  
\end{tm}

\noindent
{\em Proof.} It suffices to bound the degree $d$ of any such $X$. 
The normal bundle $N$ of $X$ in $\Q{n}$ is generated by global sections.
The Proposition of \ci{sc} is valid, on $X$, with ``ample" replaced
by ``generated by global sections;" see \ci{fu}, Example 12.1.7.
It follows that
\begin{equation}
\label{schur}
c_{n-\nu}(N) \cdot c_1(N)^{2\nu -n} \leq c_1(N)^{\nu}.
\end{equation}
By the self intersection formula for $X$
on $\Q{n}$ and   the structure of the cohomology
ring of quadrics  we have 
$c_{n-\nu}(N)=\frac{1}{2}dL^{n -\nu}$. 

\noindent
By \ci{be-so-book}, Theorem 2.3.11 we get that
$Pic(X)\simeq \zed [L]$, so that, if $X$ is not of  general type,
$K_X=eL$, with $e\leq 0$. Adjunction formula gives 
$c_1(N)=(e+n)L$.
By plugging into (\ref{schur}) we get
$$
\frac{1}{2}d^2 (e+n)^{2\nu-n} \leq (e+n)^{\nu}d,
$$
which gives, after observing that $0 \leq  -e \leq \dim X +1 < n $, that
$$
d\leq 2 (e+n)^{n-\nu}\leq 2n^{n-\nu}.
$$
\blacksquare

\section{$3$-folds on  a hypersurface of fixed degree}
\label{chTHFD}
In this section we generalize to the case of
$\Q{5}$ the main result of $\S 3$ of \ci{boss1},
which deals with bounds associated with
  nonsingular $3$-folds contained in a hypersurface of 
$\pn{5}$. For the analogous result on $\Q{4}$
see \ci{a-s}, Proposition $6.7$. However in both 
of the above references
the result is proved under the assumption that ``$d$ 
is big enough" with
respect to the degree of the hypersurface. 
Of course this
assumption is not a real restriction, since the residual cases
are automatically taken care
of by the fact that having a bounded degree bounds everything. 
However, it seems
convenient to prove our statements without restrictions.

The importance of this bound is more ore less theoretical: 
it can be used to assert
the finiteness of 
special families of $3$-folds in $\Q{5}$. One should not expect
to make an effective use of them  and get sharp results. 
The paper
\ci{e-p}, which deals with surfaces in $\pn{4}$, is the 
original source
of the main ideas used in \ci{boss1}, in $\S6$ of  \ci{a-s}, 
and in this section.
The theoretical bound given there, for the degree of
nonsingular surfaces not of general type in $\pn{4}$, 
is not an effective one.
In the paper \ci{b-f} an effective bound, $d\leq 105$, 
is proved using
initial ideals.

\bigskip
Let $X$ be a $3$-fold of degree $d$ in $\Q{5}$ contained in 
an  integral hypersurface 
$V\in |\odixl{\Q{5}}{\s}|$, $S$ a general hyperplane 
section of $X$, $C$
a general hyperplane section of $S$ and let $g$ be its genus.
As a convention, when we write something like
``$+$ l.t. in $\sqrt{d}$," we mean that the coefficients of 
the lower terms depend
only on $\s$.

\begin{tm}
\label{polybound}
Let 
$X\subseteq V\subseteq\Q{5}$ 
be as above. There is a degree eight polynomial
$P_{\s}(\sqrt d)$, depending only on
$\s$ and with positive leading coefficient, such that
$$
-\chi(\odix{X})\geq P_{\s}(\sqrt d).
$$

\end{tm}

\noindent
{\em Proof.}
\noindent
Look at the following three exact sequences.
$$
0\to \odixl{\pn{6}}{t-2} \to \odixl{\pn{6}}{t}\to 
\odixl{\Q{5}}{t} \to 0,
$$
$$
0\to \odixl{\Q{5}}{t-\s} \to \odixl{\Q{5}}{t}\to 
\odixl{V}{t} \to 0,
$$
$$
0\to {\cal I}_{X,V}(t) \to \odixl{V}{t}\to \odixl{X}{t} 
\to 0.
$$
One can use the first one to compute
$\chi (\odixl{\Q{5}}{t})$,
 the second one to compute
\begin{eqnarray*}
\chi (\odixl{V}{t}) &= & \frac{1}{12}\s t^4 + (-\frac{1}{6}\s^2+
 \frac{5}{6}\s)t^3
+(\frac{1}{6}\s^3- \frac{5}{4}\s^2+3\s)t^2     \\
& & +(-\frac{1}{12}\s^4+ \frac{5}{6}\s^3- 3\s^2 + \frac{55}{12}\s)t  
 +\frac{1}{60}\s^5 - \frac{5}{24}\s^4 + \s^3 - \frac{55}{24}\s^2 
+ \frac{149}{60}\s,
\end{eqnarray*}
finally we use (\ref{chioxt}),  
$\mu:=\mu_{\s}=$
$\frac{1}{2}d^2+ \s(\s-3)d -2\s(g-1)$ (cf. the notation fixed 
before Lemma \ref{epas}),
 and the third short exact sequence to compute
\begin{eqnarray*}
\chi ({\cal I}_{X,V}(t)) &= & 
\frac{1}{12}\s t^4 + \frac{1}{6}[(5-\s)\s -d]t^3+   \\
& &[\frac{1}{6}\s^3 -\frac{5}{4}{\s}^2 + 3\s  +
 \frac{1}{4\s}(\frac{d^2}{2} + d\s (\s-3) -\mu ) -\frac{d}{2}]t^2 
  \\
& &  +[ -\frac{1}{12}\s^4 + \frac{5}{6}{\s}^3-3{\s}^2 + \frac{55}{12}\s 
-\frac{d}{3} \\   
& & + \frac{1}{4\s}(\frac{d^2}{2} + d\s (\s-3) -\mu) -
\chi (\odix{S})]t   \\
& & + \frac{1}{60}\s^5 - \frac{5}{24}\s^4 + \s^3 -\frac{55}{24}\s^2
+ \frac{149}{60}\s -\chi (\odix{X})  \\
&=:& Q(t)-\chi (\odix{X}).
\end{eqnarray*}
It follows that
$$
-\chi (\odix{X})=\chi ({\cal I}_{X,V}(t)) - Q(t).
$$
Define
$$
t_1:=\min\{ t\in {\Bbb N}|\ \ \delta:=2\s t-d>0\ \  
{\rm and}\ \ \frac{\delta^2}{2} -\mu -\delta \s (\s-3)>0\}; 
$$
by \ci{a-s}, page 89:
$$
\frac{d}{2\s}\leq t_1 \leq \frac{d}{2\s} + \frac{\sqrt{2d}}{2} + \s
$$
By plugging $t_1$ in what above we get, using the above 
inequalities for $t_1$
and Lemma \ref{epas}:
$$
-\chi (\odix{X})=\chi ({\cal I}_{X,V}(t_1)) - Q(t_1)\geq
\chi ({\cal I}_{X,V}(t_1))-\frac{1}{64{\s^3}}d^4 + 
\frac{1}{2\s}\chi (\odix{S})d +\  {\rm l.t.\ in}\ \sqrt{d};
$$
by \ci{a-s}, pages $88-89$, we also know that
$$
\chi (\odix{S})\geq \frac{1}{24{\s}^2}d^3 +\ {\rm l.t.\ in}\sqrt{d}
$$
so that
$$
-\chi (\odix{X})\geq \chi ({\cal I}_{X,V}(t_1))  +
\frac{1}{192{\s}^3}d^4 + \ {\rm l.t.\ in}\sqrt{d}.
$$
To conclude it is enough to bound conveniently from below
$\chi ({\cal I}_{X,V}(t_1))=h^0-h^1+h^2-h^3+ h^4$. 
This, in turn, can be accomplished by bounding $h^1$ and $h^3$ from above.
This is the content of the following
technical lemmata. 
\blacksquare

First we fix some notation.
By taking general hyperplane sections we obtain the following diagram:

$$
\begin{array}{rllllllll}
 \hspace{1in}\Q{3} & \subseteq & \ldots & \subseteq & \Q{n+1} & 
\subseteq & 
\Q{n+2} & \  & \  \\
  \cup \ \  & \ & \ & \ & \cup & \ & \cup & \ & \   \\
\tilde{V}^2\to V^2  & \subseteq & \ldots & \subseteq & V^n &
\subseteq & V^{n+1} & = & V  \\
 \cup \ \  & \ & \ & \ & \cup & \ & \cup & \ & \   \\
X^1 & \subseteq & \ldots & \subseteq & X^{n-1} & \subseteq & X^n &
= & X
\end{array}
$$
where $\tilde{V}^2$ is the normalization of $V^2$.

\medskip
The following lemma is the analogue of
\ci{boss1}, Lemma $3.3$. It is proved in the same way
using \ci{a-s} lemmata $6.10$, $6.11$ and $6.12$ instead of the
lemmata from \ci{e-p} quoted in  \ci{boss1}.

\begin{lm}
\label{uno}
Let $X=X^n$ be a degree $d$ 
 nonsingular  $n$-dimensional subvariety of $\Q{n+2}$, $n\geq 2$.
Assume that $X$ is contained in an  integral hypersurface 
$V=V^{n+1}$  $\in |\odixl{\Q{n+2}}{\s}|$.
Define $t_1$ as above.
Then there are constants $A_1$, $A_2$, depending only on $\s$, 
such that
$$
\sum_{\nu=t_1}^{\infty} h^1({\cal I}_{X^1,\tilde{V}^2}(\nu))
\leq A_1 \sqrt{d^3} + \ \rm{l.t.\ in }\  \sqrt{d},
$$
and
$$
\sum_{\nu=0}^{t_1-1} h^1({\cal I}_{X^1,\tilde{V}^2}(\nu))
\leq A_2 \sqrt{d^5} + \ \rm{l.t.\ in }\  \sqrt{d}.
$$
\end{lm}

\medskip
The next lemma  merely generalizes Lemma 3.4 of \ci{boss1}. 
It should be noted that
their proof of it has a flaw since their argument does not 
work in the case
$i=n+1$. However that case is not needed for our (and their) 
purposes. In any case
we  easily prove a bound also in that case.

\begin{lm}
\label{due}
Let things be as in the previous lemma. Then
$$
h^0({\cal I}_{X,V}(t_1))\leq B_0\sqrt{d^{2n-1}} +\ {\rm  l.t. \ in }\ 
\sqrt{d};
$$ 
for 
$i=1,$
$$
h^1({\cal I}_{X,V}(t_1))\leq B_1\sqrt{d^{2n+1}} +\ {\rm  l.t. \ in }\ 
\sqrt{d},
$$
for $i=n-1,\ n$
$$
h^i({\cal I}_{X,V}(t_1))\leq B_i\sqrt{d^{2i+1}} +\ {\rm  l.t. \ in }\ 
\sqrt{d},
$$ 
 and for
$i=n+1$
$$
h^{n+1}({\cal I}_{X,V}(t_1))\leq B_{n+1} d^{n+1} +\ {\rm  l.t. \ in }\ 
\sqrt{d},
$$
where the $B_i$'s are positive constants depending only on $\s$.
\end{lm}

\noindent
{\em Proof.}
By looking at  the following sequences
\begin{equation}
\label{alce}
0\to
{\cal I}_{X^{i},V^{i+1}}(k-1)
\to
{\cal I}_{X^{i},V^{i+1}}(k) 
\to
{\cal I}_{X^{i-1},V^{i}}(k)
\to
0,
\end{equation}
we deduce
\begin{eqnarray*}
h^0({\cal I}_{X,V}(t_1)) & \leq & h^0({\cal I}_{X^{n-1},V^{n}}(t_1)) +
h^0({\cal I}_{X^{n},V^{n+1}}(t_1-1))  \\
 &\leq & \sum_{k=1}^{t_1} h^0({\cal I}_{X^{n-1},V^{n}}(k))  \\
& \vdots &   \\
 & \leq & \sum_{k=1}^{t_1} \dots \sum_{k=1}^{t_1}  
h^0({\cal I}_{X^{1},V^{2}}(k)) \qquad \quad (n-1)\  
\rm{summands}           \\
& \leq  & t_1^{n-2} \sum_{k=1}^{t_1}  h^0({\cal I}_{X^{1},V^{2}}(k))
\ \leq\  t_1^{n-2} \sum_{k=1}^{t_1}  
h^0({\cal I}_{X^{1},\tilde{V}^{2}}(k))  \\
& \leq & t_1^{n-2}(A_0\sqrt{d^3} + {\rm l.t \ in }\  \sqrt{d}),
\end{eqnarray*}
where the last inequality can be found in  \ci{a-s}, 
Lemma  6.15, and  $A_0$
depends only on $\s$.
Since $t_1\leq \frac{1}{2\s}d + \frac{\sqrt{2d}}{2} + \s$,
we have bounded 
$h^0({\cal I}_{X,V} (t_1))$ 
as wanted.

\smallskip
\noindent
To bound $h^1$ we argue as above. 

\noindent
For
$h^1({\cal I}_{X,V} (t_1))\leq$ 
$\sum_{k=1}^{t_1}h^1({\cal I}_{X^{n-1},V^n} (k))+$
$h^1({\cal I}_{X^n,V^{n+1}})$, but this last dimension is
zero as one can check by looking  at the  long
cohomology sequences associated with the following two exact sequences:
\begin{equation}
\label{aaa}
0\to {\cal I}_{X^n,\Q{n+2}}
\to
\odix{\Q{n+2}}
\to
\odix{X^n}
\to 0,
\end{equation}
\begin{equation}
\label{aaaa}
0 \to
\odixl{\Q{n+2}}{-\s}
\to
{\cal I}_{X^n,\Q{n+2}}
\to
{\cal I}_{X^n,V^{n+1}}
\to 0.
\end{equation}
An easy induction argument, already seen before, using
(\ref{alce}) allows us to infer that
$$
h^1({\cal I}_{X,V}(t_1))\leq t_1^{n-2}
\sum_{k=1}^{t_1}h^1({\cal I}_{X^1,V^2}(k)).
$$
To obtain the desired bound on $h^1$ it is enough
to prove that:
$$
\sum_{k=1}^{t_1}h^1({\cal I}_{X^1,V^2}(k)) \leq
F\sqrt{d^5} +\ \rm{l.t. \ in} \sqrt{d},
$$
where $F$ is a constant depending only on $\s$.

\noindent
This can be proved as follows. The idea is to couple
the previous lemma with the  cohomology sequences associated with
the following exact sequences:
$$
0\to
{\cal I}_{X^1,V^2}(k)\to
{\cal I}_{X^1,\tilde{V}^2}(k)\to 
 Q(k)\to 0.
$$
\noindent
Clearly we have
$
h^1({\cal I}_{X^1,V^2}(k))\leq 
h^1({\cal I}_{X^1,\tilde{V}^2}(k))+ h^0(Q(k))
$, $\forall k$, so that, in view of Lemma
\ref{uno},
it is enough to prove the following 

\medskip
\noindent
CLAIM. $h^0(Q(k))\leq D(k+1)$, $\forall k\geq 0$, {\em where
$D$ is a positive constant depending only on $\s$. In particular
$\sum_{k=0}^{t_1}h^0(Q(k))\leq (1/2)Dt_1^2+ $ \rm{l.t.} in $t_1$.}

\medskip
\noindent
{\em Proof of the claim.} $Q$ is the structural sheaf of
  the non-normal locus of $V^2$ twisted
 by the ideal sheaf of $X^1$. By taking a general hyperplane section we 
get the following exact sequences

$$
0 \to
Q(k-1)
\to
Q(k)
\to
Q_{\Gamma}(k)
\to
0,
$$
where 
$Q_{\Gamma}(k)\simeq Q_{\Gamma}$ is the  structural
 sheaf of the singular locus
of $\Gamma$, a general hyperplane section of
$V^2$.
As usual
$h^0(Q(k))\leq h^0(Q(k-1)) + length(Q_{\Gamma})$, so that
$h^0(Q(k))\leq length(Q_{\Gamma})(k+1) + h^0(Q(-1))$.
This length is bounded from above by a function of $\s$ only 
(that is,
an irreducible curve of degree $2\s$ on a two dimensional quadric
cannot have too many singularities).
As to $ h^0(Q(-1))$ one has to exercise caution since
the non-normal locus may  be  non reduced. However by  looking
at the cohomology
sequences associated with
(\ref{aaa})
 we see that
$h^1({\cal I}_{X^1,V^2}(-1))=0$, so that $ h^0(Q(-1))=0$, 
and the claim is proved.

\smallskip
\noindent
Now we prove the bound for $i=n$.

\noindent
We  start by remarking that for
$i=n-1,$ $n,$ $n+1$ and all $k\geq d-1$.
$$
h^i({\cal I}_{X,V}(k))=0.
$$
 If
$X^2$ is non degenerate then 
${\cal I}_{X^2,\pn{5}}$ is 
$(d-2)$-regular (cf. \cite{la2}) in the sense of Castelnuovo-Mumford.
By looking at the following sequences:
$$
0\to
{\cal I}_{X^2,\pn{5}}(-2+k)
\to
{\cal I}_{X^2,\pn{5}}(k)
\to
{\cal I}_{X^2,\Q{4}}(k)
\to 0,
$$
we deduce the vanishings, for $n=2$.
An easy inductive argument (cf. \ci{boss1}, page 326) 
gives the desired vanishings. 

\noindent
If $X^2
\subseteq \pn{5}$ were degenerate, then either $X^2=\pn{2}$, 
$X^2$ would be
 a hypersurface
of $\pn{3}$ or it would be
 a nondegenerate surface in $\pn{4}$. In any of these cases we apply
the bound for the regularity of the ideal sheaves
in \ci{la2} to obtain vanishings for the higher cohomology
of ${\cal I}_{X^2,{\Bbb P}^{4}}$ which is easy to lift to the 
desired vanishings
for $X^2$. Again the inductive procedure allows us to conclude. 

\noindent
We have the following chain of inequalities:

\begin{eqnarray*}
h^n({\cal I}_{X^n,V^{n+1}}(t_1)) & \leq &
\sum_{k>t_1}h^{n-1}({\cal I}_{X^{n-1},V^{n+1}}(k))
\leq \sum_{k=1}^{d-4}
h^{n-1}({\cal I}_{X^{n-1},V^{n+1}}(k))  \\
 & \leq & \ldots \leq  
\sum_{1}^{d-4} \ldots \  \sum_{1}^{d-4} 
\sum_{k=1}^{d-4} h^1({\cal I}_{X^1,V^2}(k))   \\
 & \leq & (d-4)^{n-2} 
\sum_{k=1}^{d-4} h^1({\cal I}_{X^1,V^2}(k))    \\
& \leq & B_n\sqrt{d^{2n+1}} + \rm{l.t.\ in} \sqrt{d},
\end{eqnarray*}
where $B_n$ depends only on $\s$ and
 the last inequality follows from Lemma \ref{uno}. 
\smallskip
\noindent
The case $i=n-1$ is analogous.

\smallskip
\noindent
Finally the bound for the case $i=n+1$ can be obtained as  in the case 
$i=n$ except for the fact that we end up having to bound
$h^2( {\cal I}_{X^1,V^2}(k))$ for $k=1,\ldots,$ $d-4$, and not $h^1$:
$$
h^{n+1}({\cal I}_{X,V}(t_1))\leq
d^{n-2}(\sum_{k=1}^{d-4}h^2({\cal I}_{X,V}(t_1)).
$$
To bound this summand we look at the exact sequences:
$$
0\to
{\cal I}_{X^1,V^2}(k)
\to
\odixl{V^2}{k}
\to
\odixl{X^1}{k}
\to
0,
$$
and deduce
$$
h^2({\cal I}_{X^1,V^2}(k))\leq h^1(\odixl{X^1}{k}) + h^2(\odixl{V^2}{k})
=h^0(\omega_{X^1}(-k))+h^0(\odixl{X^1}{-3+\s -k})
$$
where the last equality stems from Serre Duality.
We are thus left with bounding the two $h^0$ above.
The first one can be bounded using Proposition \ref{easybound}
on $ h^0(\omega_{X^1})=g(X^1)$: the worst upper bound is
of the form
$(1/4\s)d^2+$ l.t. in $d$.
As to the second $h^0$ its worst upper bound is of the form
$(1/2)\s^2$. Adding up for $k=1,\ldots,$ $d-4$ we get
that the worst upper bound is
$(1/4\s)d^3 +$ l.t. in $d$.
\blacksquare

\bigskip
The following generalizes \ci{boss1}, Corollary $3.1$. 

\begin{cor}
\label{finiteonhyp}
Let $\s$ be any positive integer. There are only finitely
 many components
 of the Hilbert scheme of
$\Q{5}$ corresponding to nonsingular $3$-folds in $\Q{5}$ 
which are not of general type and
are contained in some hypersurface of degree $\s$.
\end{cor}

\noindent{\em Proof.} It is enough to bound from 
above the degree of such $3$-folds.
Since $\omega_X(-1)$ does not have sections
$h^0(\omega_X)\leq h^0(\omega_S)$, where $S$ is any nonsingular 
hyperplane section
of $X$. By the generalized Castelnuovo-type bounds of Harris 
(cf. \ci{jh})
we have
$$
h^0(\omega_S)\leq Ad^3 + \rm{l.t.\ in }\ d,
$$ 
where $A$ is some  constant;
the Lefschetz hyperplane theorem, coupled with Proposition
\ref{easybound}, ensures
that
$$
h^1(\odix{X})= h^1(\odix{S})\leq h^1(\odix{C})\leq \frac{1}{4\s}d^2 +
\rm{l.t.\ in }\ d. 
$$
It follows that
$$
h^0(\omega_S) \geq h^0(\omega_X)= 1+ h^2(\odix{X}) -h^1(\odix{X}) -
\chi(\odix{X})\geq \frac{1}{192\s^3}d^4 + \rm{l.t.\ in }\ \sqrt{d}.
$$ 
Comparing the two inequalities for $h^0(\omega_S)$, we conclude 
that $d$ is bounded.
\blacksquare

\section{Finiteness for $3$-folds not of general type in $\Q{5}$}
\label{boundednessonq5}

\begin{pr}
\label{generalizedhodge}
Let $X$ be a nonsingular $3$-fold in $\Q{5}$ and $k$ a positive 
integer. Then
\begin{equation}
\label{genhodge}
\chi (\odix{S})\leq \frac{2}{3} \frac{(g-1)^2}{d} - \frac{1}{24}d^2 + 
\frac{5}{12}d.
\end{equation}
If $X$ is not of general type then
\begin{equation}
\label{Xnotgentype-chi}
-\chi (\odix{X}) \leq \chi (\odix{S}) + \frac{1}{2}d^2 -2d +2;
\end{equation}
if $d>2k^2$ and $X$ is not of general type and not contained
in any hypersurface of degree strictly less than
$2\cdot k$, then
\begin{equation}
\label{Xnotgentype-chik}
-\chi (\odix{X}) \leq \chi (\odix{S}) + \frac{1}{k}d^2 + (k-4)d +2
\end{equation}
\end{pr}

\noindent
{\em Proof.}
The first inequality stems from the Generalized Hodge Index Theorem
contained in  \ci{boss1}:
$$
d(K_X^2L)\leq (K_XL^2)^2,
$$
we make explicit the left hand side using
(\ref{K2L}) and the right hand side using
(\ref{KL2}).

\noindent
For the second one we look at
$$
0 \to K_X(-1) \to K_X \to K_S(-1) \to 0.
$$
Since $X$ is not of general type
$h^0(K_X(-1))=0$, otherwise $K_X$ would be big, i.e. a 
$|mK_X|$ would define a birational map. It follows that
$h^3(\odix{X})=$ $h^0(K_X)$ $\leq$ $ h^0(K_S(-1))$ 
 $\leq$ $ h^0(K_S)=$ $ h^2(\odix{S})$.

\noindent
We thus have
$$
-\chi (\odix{X})\leq  h^1(\odix{X}) + h^3(\odix{X})
\leq \chi (\odix{S}) + 2h^1(\odix{X}),
$$
where we have used Lefschetz Theorem on Hyperplane Sections
to ensure that $h^1(\odix{X})= h^1(\odix{S})$. $h^1(\odixl{S}{-1})=0$
by Kodaira Vanishing, so that
 $h^1(\odix{X})= h^1(\odix{S})\leq g$.

\noindent
If $C$ were contained in a $\pn{3}$ we would use
Proposition \ref{easybound} with $k=1$ to conclude.
If $C$ were not in any surface of degree  strictly less than
$2\cdot 2$ we would use Proposition \ref{boundasep} with $k=2$.

\noindent
The third inequality is proved exactly as the second one by using
 Proposition \ref{roth} to ensure that a general curve section $C$
is not in any surface of the corresponding $\Q{3}$ of degree 
strictly less than
$2\cdot k$, and Proposition \ref{boundasep} to bound the genus 
$g$ from above. 
\blacksquare 

\begin{pr}
\label{s3>=0forN}
Let $X$ be a nonsingular $3$-fold in $\Q{5}$.
Then
$$
60 \chi (\odix{S}) \geq \frac{3}{2}d^2 - 12d + (d-48)(g-1) + 24
\chi (\odix{X}).
$$
\end{pr}

\noindent
{\em Proof.}
Denote by $s_i$ and $n_i$ the Segre and Chern classes respectively
 of the normal bundle 
$N$ of $X$ in $\Q{5}$. Since $N$ is generated by global sections
we have $s_3\geq 0$. Since $s_3=n_1^3- 2n_1n_2$ we get

\noindent
$0\leq (K_X+5L)^3 - 2(K_X+5L)\frac{1}{2}dL^2=$
$K^3 + 15K_X^2L + 75K_XL^2 + 125d - d(K_X+5L)L^2.$

\noindent
We conclude by
(\ref{K3}), (\ref{K2L}) and (\ref{KL2}).
\blacksquare

\begin{tm}
\label{sigh}
Let $n=4,$ $5$ or $n\geq 7$.
There are only finitely many components of the Hilbert scheme
of $\Q{n}$ corresponding to nonsingular $(n-2)$-folds not 
of general type.
\end{tm}

\noindent
{\em Proof.}
By \ci{a-s}, \S6 and Theorem \ref{boundednessonqn} it 
is enough to consider
the case  $n=5$.
It is enough to bound from above  the degree $d$ of such $3$-folds.

\noindent
 Fix a positive integer $k$ and let $d$ be a positive integer 
such that
 $d >2k^2.$
Let  $X$ be  a  degree $d$ $3$-fold in $\Q{5}$
not lying on any hypersurface of $\Q{5}$ of degree 
 strictly less than $2\cdot k$;  by Proposition \ref{roth},
  a general curve section of $X$
does not lie on any surface of the corresponding $\Q{3}$ of 
degree strictly less
than
$2\cdot k$.

\noindent
 We couple Proposition \ref{s3>=0forN} and inequality
(\ref{Xnotgentype-chik}) of Proposition \ref{generalizedhodge}:
$$
84 \chi (\odix{S}) \geq (\frac{3}{2} - \frac{24}{k})d^2 -
( 24k-84)d -48 + (d-48)(g-1).
$$
In what above we plug inequality (\ref{genhodge}) of
Proposition \ref{generalizedhodge} and get:
$$
\frac{52}{d} (g-1)^2 - \frac{21}{6}d^2 +35d \geq (\frac{3}{2}- 
\frac{24}{k})d^2 - ( 24k -84)d -48 + (d-48)(g-1).
$$
A simple manipulation gives
\begin{equation}
\label{almost}
(g-1)[ \frac{52}{d}(g-1) -d +48] + (\frac{24}{k}-5)d^2 +
(24k - 49)d +48 \geq0.
\end{equation}
Let us now first assume that $g>0$. The aim is to choose $k$ such that
the coefficients
$\alpha:=(\frac{52}{d}(g-1) -d +48)$ and 
$\beta:=(\frac{24}{k} -5)$ of $(g-1)$ and $d^2$, respectively, 
are negative.
Once they are negative, since $k$ is fixed, the  inequality
(\ref{almost})
will force $d$ to be bounded from above.
 By Proposition \ref{epas} we get
\begin{eqnarray*}
\frac{52}{d}(g-1) - d + 48 & \leq & 
\frac{52}{d}[\frac{d^2}{2k} + \frac{1}{2}(k-4)d] -d +48  \\
& = & (\frac{26}{k} -1)d + 26k -56.
\end{eqnarray*}
Let $k=27$; then $\alpha$ is negative for $d\gg 0$.
For the same value of $k$, $\beta$ is negative as well.
By what above we infer that $d$ is bounded from above if $g>0$ and
$X$ is not in an hypersurface of degree strictly less than
$2\cdot 27$. We apply Corollary \ref{finiteonhyp}
to see that $d$ is bounded from above for $3$-folds $X$, 
not of general type,
contained in  hypersurfaces of degrees less or equal to $2\cdot 27$.
This proves the theorem in the case $g>0$.

\noindent
Assume that $g=0$. Then, by (\ref{almost}):
$$
\frac{52}{d} + d -48 + (\frac{24}{k} -5)d^2 + (24k -49)d +48 \geq 0.
$$
 We argue as above with 
$k=5$.
\blacksquare

\end{document}